\input amstex
\input amsppt.sty
\magnification=\magstep1
\hsize=32truecc
\vsize=22.5truecm
\baselineskip=16truept
\NoBlackBoxes
\TagsOnRight \pageno=1 \nologo
\def\Z{\Bbb Z}

\def\l{\left}
\def\r{\right}
\def\bg{\bigg}
\def\({\bg(}
\def\[{\bg\lfloor}
\def\){\bg)}
\def\]{\bg\rfloor}
\def\t{\text}
\def\f{\frac}

\def\sm{\setminus}

\def\eq{\equiv}

\def\ls{\leqslant}
\def\gs{\geqslant}
\def\mo{\roman{mod}}
\def\ord{\roman{ord}}

\def\da{\delta}

\def\Proof{\noindent{\it Proof}}

\def\Remark{\medskip\noindent{\it  Remark}}

\def\Ack{\medskip\noindent {\bf Acknowledgments}}
\hbox {Acta Arith. 190 (2019), no.\,3, 209--220.}
\bigskip
\topmatter
\title Consecutive primes and Legendre symbols\endtitle
\author Hao Pan and Zhi-Wei Sun\endauthor
\leftheadtext{Hao Pan and Zhi-Wei Sun}
 \rightheadtext{Consecutive primes and Legendre symbols}
\address
 (Hao Pan) School of Applied Mathematics, Nanjing University of Finance and Economics,
 Nanjing 210046, People's Republic of China
\endaddress
\email haopan79\@zoho.com\endemail
\address
(Zhi-Wei Sun, corresponding author) Department of Mathematics, Nanjing University,
 Nanjing 210093, People's Republic of China
\endaddress
\email zwsun\@nju.edu.cn\endemail

\abstract Let $m$ be any positive integer and let $\da_1,\da_2\in\{1,-1\}$. We show that for some constant $C_m>0$ there are infinitely many integers $n>1$
with $p_{n+m}-p_n\ls C_m$ such that
$$\l(\f{p_{n+i}}{p_{n+j}}\r)=\da_1\ \quad\t{and}\ \quad\l(\f{p_{n+j}}{p_{n+i}}\r)=\da_2$$
for all $0\ls i<j\ls m$, where $p_k$ denotes the $k$-th prime, and $(\f {\cdot}p)$ denotes the Legendre symbol for any odd prime $p$.
We also prove that under the Generalized Riemann Hypothesis there are infinitely many positive integers $n$ such that $p_{n+i}$
is a primitive root modulo $p_{n+j}$ for any distinct $i$ and $j$ among $0,1,\ldots,m$.
\endabstract
\thanks 2010 {\it Mathematics Subject Classification}. Primary 11A41, 11N05; Secondary 11A07, 11A15.
\newline\indent The first author is supported by the National Natural Science
Foundation of China (grant 11671197), and the second author is supported by the National Natural Science
Foundation of China (grant 11571162) and the NSFC-RFBR Cooperation and Exchange Program (grant no. 11811530072).
\endthanks
\endtopmatter
\document

\heading{1. Introduction}\endheading

For $n\in\Z^+=\{1,2,3,\ldots\}$ let $p_n$ denote the $n$-th prime.
The famous twin prime conjecture asserts that $p_{n+1}-p_n=2$ for infinitely many $n\in\Z^+$.
Although this remains open, recently Y. Zhang [Z] was able to prove that
$$\liminf_{n\to\infty}(p_{n+1}-p_n)\ls 7\times10^7.$$
The upper bound $7\times 10^7$ was later reduced to 4680 by the {\tt Polymath} team [Po] led by T. Tao, and 600 by J. Maynard [M], and 246 by the {\tt Polymath} team [Po].
Moreover, J. Maynard [M], as well as T. Tao, established the following deep result.

\proclaim{Theorem 1.1 {\rm (Maynard-Tao)}} For any positive integer $m$, we have
$$\liminf_{n\to\infty}(p_{n+m}-p_n)\ls Cm^3e^{4m},$$
where $C>0$ is an absolutely constant.
\endproclaim

Earlier than this work, in 2000 D.K.L. Shiu [S] proved the following nice theorem.

\proclaim{Theorem 1.2 {\rm (Shiu)}} Let $a\in\Z$ and $q\in\Z^+$ be relatively prime. Then, for any $m\in\Z^+$ there is a positive integer $n$ such that
$$p_n\eq p_{n+1}\eq\cdots\eq p_{n+m}\eq a\pmod q.$$
\endproclaim

This was recently re-deduced in [BFTB] via the Maynard-Tao method.

In this paper we mainly establish the following new result on consecutive primes and Legendre symbols.

\proclaim{Theorem 1.3} Let $m$ be any positive integer and let $\delta_1,\da_2\in\{1,-1\}$. For some constant $C_m>0$ depending only on $m$,
there are infinitely many integers $n>1$ with $p_{n+m}-p_n\ls C_m$ such that for any $0\ls i<j\ls m$ we have
$$\l(\f{p_{n+i}}{p_{n+j}}\r)=\da_1\ \ \t{and}\ \ \l(\f{p_{n+j}}{p_{n+i}}\r)=\da_{2}.\tag1.1$$
\endproclaim
\Remark\ 1.4. (a) Instead of (1.1) in Theorem 1.3, actually we may require both (1.1) and the following property:
$$p_{ij}\|(p_{n+i}-p_{n+j})\quad\t{for some prime}\ p_{ij}>2m+1.\tag1.2$$
(As usual, for a prime $p$ and an integer $a$, by $p\|a$ we mean $p\mid a$ but $p^2\nmid a$.)

(b) We conjecture the following extension of Theorem 1.3: For any $m\in\Z^+$, $\da\in\{1,-1\}$ and  $\da_{ij}\in\{1,-1\}$ with $0\ls i<j\ls m$,
there are infinitely many integers $n>1$ such that
$$\l(\f{p_{n+i}}{p_{n+j}}\r)=\da_{ij}=\da\l(\f{p_{n+j}}{p_{n+i}}\r)$$
for all $0\ls i<j\ls m$.
\medskip

{\it Example} 1.5. The smallest integer $n>1$ with
$$\l(\f{p_{n+i}}{p_{n+j}}\r)=1\quad\t{for all}\ i,j=0,\ldots,6\ \t{with}\ i\not=j$$
is $178633$, and a list of the first 200 such values of $n$ is available from [Su, A243901].
The seven consecutive primes
$p_{178633},\ p_{178634},\ \ldots,\ p_{178639}$ have the concrete values
$$2434589,\ 2434609,\ 2434613,\ 2434657,\ 2434669,\ 2434673,\ 2434681$$
respectively.
\medskip

{\it Example} 1.6. The smallest integer $n>1$ with
$$\l(\f{p_{n+i}}{p_{n+j}}\r)=-1\quad\t{for all}\ i,j=0,\ldots,5\ \t{with}\ i\not=j$$
is $2066981$, and the six consecutive primes
$p_{2066981},\ p_{2066982},\ \ldots,\ p_{2066986}$ have the concrete values
$$33611561,\ 33611573,\ 33611603,\ 33611621,\ 33611629,\ 33611653$$
respectively.
\medskip

{\it Example} 1.7. The smallest integer $n>1$ with
$$-\l(\f{p_{n+i}}{p_{n+j}}\r)=1=\l(\f{p_{n+j}}{p_{n+i}}\r)\quad\t{for all}\ 0\ls i<j\ls 6$$
is $7455790$, and the seven consecutive primes
$p_{7455790},\ p_{7455791},\ \ldots,\ p_{7455796}$ have the concrete values
$$131449631,\, 131449639,\, 131449679,\, 131449691,\, 131449727,\, 131449739,\, 131449751$$
respectively.
\medskip

{\it Example} 1.8. The smallest integer $n>1$ with
$$\l(\f{p_{n+i}}{p_{n+j}}\r)=1=-\l(\f{p_{n+j}}{p_{n+i}}\r)\quad\t{for all}\ 0\ls i<j\ls 5$$
is $59753753$, and the six consecutive primes
$p_{59753753},\ p_{59753754},\ \ldots,\ p_{59753758}$ have the concrete values
$$1185350899,\ 1185350939,\ 1185350983,\ 1185351031,\ 1185351059,\ 1185351091$$
respectively.
\medskip

Actually Theorem 1.3 is motivated by the following conjecture of the second author.

\proclaim{Conjecture 1.9 {\rm (Sun [Su, A243837])}} For any positive integer $m$, there are infinitely many $n\in\Z^+$ such that for any distinct $i$ and $j$ among $0,1,\ldots,m$
the prime $p_{n+i}$ is a primitive root modulo $p_{n+j}$.
\endproclaim

{\it Example} 1.10. The least $n\in\Z^+$ with $p_{n+i}$ a primitive root modulo $p_{n+j}$ for any distinct $i$ and $j$ among $0,1,2,3$ is 8560,
and a list of the first 50 such values of $n$ is available from [Su, A243839].
Note that
$$p_{8560}=88259,\ p_{8561}=88261\ \t{and}\ \ p_{8562}=88289.$$

Our second result confirms Conjecture 1.9 under the Generalized Riemann Hypothesis.

\proclaim{Theorem 1.11} Let $m$ be any positive integer. Assuming the GRH (Generalized Riemann Hypothesis). for some constant $C_m>0$ depending only on $m$, there are infinitely many integers $n>1$ with $p_{n+m}-p_n\ls C_m$, such that the prime $p_{n+i}$ is a primitive root modulo $p_{n+j}$ for any distinct $i,j\in\{0,1,\ldots,m\}$.
\endproclaim

We will prove Theorem 1.3 in the next section with the help of the Maynard-Tao work, and show Theorem 1.11 in Section 3 by combining our method with
a recent result of P. Pollack [P] motivated by the Maynard-Tao work on bounded gaps of primes and Artin's conjecture on primitive roots modulo primes.

Throughout this paper, $p$ always represents a prime. For two integers $a$ and $b$, their greatest common divisor is denoted by $\gcd(a,b)$.

\heading{2. Proof of Theorem 1.3}\endheading

Let $h_1,h_2,\ldots,h_k$ be distinct positive integers. If $\bigcup_{j=1}^kh_i(\mo\ p)\not=\Z$ for any prime $p$
(where $a (\mo\ p)$ denotes the residue class $a+p\Z$), then we call $\{h_i:\ i=1,\ldots,k\}$ an {\it admissible set}.
Hardy and Littlewood conjectured that if $\Cal H=\{h_i:\ i=1,\ldots,k\}$ is admissible then there are infinitely many $n\in\Z^+$
such that $n+h_1,n+h_2,\ldots,n+h_k$ are all prime. We need the following result in this direction.

\proclaim{Lemma 2.1 {\rm (Maynard-Tao)}} Let $m$ be any positive integer. Then there is an integer $k>m$ depending only on $m$ such that
if $\Cal H=\{h_i:\ i=1,\ldots,k\}$ is an admissible set of cardinality $k$ and $W=q_0\prod_{p\ls w}p\ ($with $q_0\in\Z^+)$ is relatively prime to $\prod_{i=1}^kh_i$
 with $w=\log\log\log x$ large enough, then for some integer $n\in[x,2x]$ with $W\mid n$ there are more than $m$ primes among
 $n+h_1,n+h_2,\ldots,n+h_k$.
\endproclaim

\proclaim{Lemma 2.2} Let $k>1$ be an integer. Then there is an admissible set $\Cal H=\{h_1,\ldots,h_k\}$ with $h_1=0<h_2<\ldots<h_k$
which has the following properties:

{\rm (i)} All those $h_1,h_2,\ldots,h_k$ are multiples of $K=4\prod_{p<2k}p$ with $p$ prime.

{\rm (ii)} Each $h_i-h_j$ with $1\ls i<j\ls k$ has a prime divisor $p>2k$ with $h_i\not\eq h_j\pmod{p^2}$.

{\rm (iii)} If $1\ls i<j\ls k$, $1\ls s<t\ls k$ and $\{i,j\}\not=\{s,t\}$, then no prime $p>2k$ divides
both $h_i-h_j$ and $h_s-h_t$.
\endproclaim
\Proof. Set $h_1=0$ and let $1\ls r<k$. Suppose that we have found nonnegative integers $h_1<\ldots<h_r$ divisible by $K$ such that
each $h_i-h_j$ with $1\ls i<j\ls r$ has a prime divisor $p>2k$ with $h_i\not\eq h_j\pmod{p^2}$, and that no prime $p>2k$ divides both $h_i-h_j$ and $h_s-h_t$
if $1\ls i<j\ls r$, $1\ls s<t\ls r$ and $\{i,j\}\not=\{s,t\}$. Let
$$X_r=\{p>2k:\ p\ \t{is prime and} \ p\mid h_s-h_t\quad\t{for some}\ 1\ls s<t\ls r\}.$$
As $K$ is relatively prime to $\prod_{p\in X_r}p$, for each $i=1,\ldots,r$ there is an integer $b_i$ with
$Kb_i\eq h_i\pmod{\prod_{p\in X_r}p}$. For each $p\in X_r$, as $r<k<p$ there is an integer $a_p\not\eq b_i\pmod p$ for all $i=1,\ldots,r$.
Choose distinct primes $q_1,\ldots,q_r$ which are greater than $2k$ but not in the set $X_r$. For any $i=1,\ldots,r$,
there is an integer $c_i$ with $Kc_i\eq h_i\pmod{q_i^2}$ since $K$ is relatively prime to $q_i^2$. By the Chinese Remainder Theorem,
there is an integer $b>h_r/K$ such that $b\eq a_p\pmod p$ for all $p\in X_r$, and $b\eq c_i+q_i\pmod {q_i^2}$ for all $i=1,\ldots,r$.

Set $h_{r+1}=Kb>h_r$. If $1\ls s\ls r$, then
$$h_{r+1}-h_s\eq Kb-Kc_s=K(b-c_s)\eq Kq_s\pmod{q_s^2},$$
hence $q_s>2k$ is a prime divisor of $h_{r+1}-h_s$ but $h_{r+1}\not\eq h_s\pmod{q_s^2}$.

For $s,t\in\{1,\ldots,r\}$ with $s\not = t$, clearly
$$\gcd(h_{r+1}-h_s,h_{r+1}-h_t)=\gcd(h_{r+1}-h_s,h_s-h_t).$$
Let $1\ls i < j\ls r$ and $1\ls s\ls r$. If a prime $p>2k$ divides $h_i-h_j$, then $p\in X_r$ and hence
$$h_{r+1}-h_s\eq Ka_p-Kb_s=K(a_p-b_s)\not\eq0\pmod p.$$
So $\gcd(h_{r+1}-h_s,h_i-h_j)$ has no prime divisor greater than $2k$.

In view of the above, we have constructed nonnegative integers $h_1<\ldots<h_k$ satisfying (i)-(iii) in Lemma 2.2.
Note that $\bigcup_{i=1}^k h_i(\mo\ p)\not=\Z$ if $p>k$. For each $p\ls k$, clearly $h_i\eq0\not\eq1\pmod p$ for any $i=1,\ldots,k$.
Therefore the set $\Cal H=\{h_1,h_2,\ldots,h_k\}$ is admissible. This concludes the proof. \qed

\medskip
\noindent{\it Proof of Theorem 1.3}. By Lemma 2.1, there is an integer $k=k_m>m$ depending on $m$ such that for any admissible set $\Cal H=\{h_1,\ldots,h_k\}$
of cardinality $k$ if $x$ is sufficiently large and $\prod_{i=1}^k{h_i}$ is relatively prime to $W=4\prod_{p\ls w}p$ then for some integer $n\in[x/W,2x/W]$
there are more than $m$ primes among $Wn+h_1,Wn+h_2,\ldots,Wn+h_k$, where $w=\log\log\log x$.

Let $\Cal H=\{h_1,\ldots,h_k\}$ with $h_1=0<h_2<\ldots<h_k$ be an admissible set satisfying the conditions (i)-(iii) in Lemma 2.2.
Clearly $K=4\prod_{p<2k}p\eq0\pmod 8.$ Let $x$ be sufficiently large with the interval $(h_k,w]$ containing more than $h_k-k$ primes. Note that $8\mid W$ since $w\gs2$. Our goal is to construct a new admissible set $\Cal H'$ to which we will apply Lemma 2.1 in order to complete the proof.

Let $\delta:=\da_1\da_2$. For any integer $b\eq\da\pmod K$ and each prime $p<2k$, clearly $b+h_i\eq\da+0\pmod p$ and hence $\gcd(b+h_i,p)=1$ for all $i=1,\ldots,k$.

By the property (ii) in Lemma 2.2, for any $1\ls i <j\ls k$, the number $h_i-h_j$ has a prime divisor $p_{ij}>2k$ with $h_i\not\eq h_j\pmod{p_{ij}^2}$.
Let $p>2k$ be an arbitrary prime dividing $\prod_{1\ls i<j\ls k}(h_i-h_j)$. Then there is a unique pair $\{i,j\}$ with $1\ls i<j\ls k$
such that $h_i\eq h_j\pmod p$. Note that $p\ls h_k$ and $p$ may be different from $p_{ij}$. All the $k-2<(p-3)/2$ numbers $h_i-h_s$ with $1\ls s\ls k$ and $s\not=i,j$ are relatively prime to $p$,
so there is an integer $r_p\not\eq h_i-h_s\pmod p$ for all $s=1,\ldots,k$ such that
$$\l(\f{r_p\,\delta}p\r)=\cases \da_2&\t{if}\ p=p_{ij},
\\1&\t{otherwise}.\endcases$$
So, for any integer $b\eq r_p-h_i\pmod p$, we have $b+h_s\not\eq0\pmod p$ for all $s=1,\ldots,k$.

Assume that $S=\{h_1,h_1+1,\ldots,h_k\}\sm\Cal H$ is a set $\{a_i:\ i=1,\ldots,t\}$ of cardinality $t>0$. Clearly $t\ls h_k-k+1$ and hence we may
choose $t$ distinct primes $q_1,\ldots,q_t\in (h_k,w]$.
If $b\eq -a_i\pmod {q_i}$, then $b+h_s\eq h_s-a_i\not\eq0\pmod {q_i}$ for all $s=1,\ldots,k$ since $0<|h_s-a_i|<h_k<q_i$.

Let
$$Q=\bg\{p\in(2k,w]:\ p\ \t{is prime and}\ p\nmid\prod_{1\ls i<j\ls k}(h_i-h_j)\bg\}\sm\{q_i:\ i=1,\ldots,t\}.$$
For any prime $q\in Q$, there is an integer $r_q\not\eq-h_i\pmod q$ for all $i=1,\ldots,k$ since $\Cal H$ is admissible.

By the Chinese Remainder Theorem, there is an integer $b$ satisfying the following (1)-(4).

(1) $b\eq \da=\da_1\da_2\pmod K$.

(2) $b\eq r_p-h_i\eq r_p-h_j\pmod p$ if $p>2k$ is a prime dividing $h_i-h_j$ with $1\ls i<j\ls k$.

(3) $b\eq -a_i\pmod {q_i}$ for all $i=1,\ldots,t$.

(4) $b\eq r_q\pmod q$ for all $q\in Q$.

By the above analysis, as we have ensured that $b+h_s\not\eq0\pmod p$ for each
 prime $p\ls w$, the product $\prod_{s=1}^k(b+h_s)$ is relatively prime to $W$.
As $\Cal H'=\{b+h_s:\ s=1,\ldots,k\}$ is also an admissible set of cardinality $k$,
for large $x$ there is an integer $n\in[x/W,2x/W]$ such that there are more than $m$
primes among $Wn+b+h_s\ (s=1,\ldots,k)$. For $a_i\in S$, we have
$$Wn+b+a_i\eq 0-a_i+a_i=0\pmod {q_i}$$
and hence $Wn+b+a_i$ is composite since $W>q_i$.
Therefore, there are consecutive primes $p_N,p_{N+1},\ldots,p_{N+m}$ with $p_{N+i}=Wn+b+h_{s(i)}$ for all $i=0,\ldots,m$, where
$1\ls s(0)<s(1)<\ldots<s(m)\ls k$. Note that
$$p_{N+m}-p_N=(Wn+b+h_{s(m)})-(Wn+b+h_{s(0)})=h_{s(m)}-h_{s(0)}\ls h_k.$$

For each $s=1,\ldots,k$, clearly $Wn+b+h_s\eq 0+\da+0=\da\pmod 8$ and hence
$$\l(\f{-1}{Wn+b+h_s}\r)=\da\ \ \t{and}\ \ \l(\f{2}{Wn+b+h_s}\r)=1,$$
where $(-)$ denotes the Jacobi symbol.
As $p_{N+i}=Wn+b+h_{s(i)}\eq\da\pmod 8$ for all $i=0,\ldots,m$, by the law of Quadratic Reciprocity we have
$$\l(\f{p_{n+j}}{p_{N+i}}\r)=\da\l(\f{p_{n+i}}{p_{N+j}}\r)\quad\t{for all}\ 0\ls i<j\ls m.$$

Let $0\ls i<j\ls m$.
Then
$$\align\l(\f{p_{N+i}}{p_{N+j}}\r)=&\l(\f{Wn+b+h_{s(i)}}{Wn+b+h_{s(j)}}\r)=\l(\f{h_{s(i)}-h_{s(j)}}{Wn+b+h_{s(j)}}\r)
\\=&\l(\f{-1}{Wn+b+h_{s(j)}}\r)\l(\f{h_{s(j)}-h_{s(i)}}{Wn+b+h_{s(j)}}\r)=\da\l(\f{h_{ij}^*}{Wn+b+h_{s(j)}}\r),
\endalign$$
where $h_{ij}^*$ is the odd part (i.e., the largest odd divisor) of $h_{s(j)}-h_{s(i)}$. For any prime divisor $p$ of $h_{ij}^*$, clearly $p\ls h_k\ls w$ and
$$\l(\f p{Wn+b+h_{s(j)}}\r)=\da^{(p-1)/2}\l(\f{Wn+b+h_{s(j)}}p\r)=\da^{(p-1)/2}\l(\f{b+h_{s(j)}}p\r).$$
If $p<2k$, then $p\mid K$, hence $b+h_j\eq\da+0\pmod p$ and thus
$$\l(\f p{Wn+b+h_{s(j)}}\r)=\da^{(p-1)/2}\l(\f{b+h_{s(j)}}{p}\r)=\da^{(p-1)/2}\l(\f{\da}p\r)=1.$$
If $p>2k$, then by the choice of $b$ we have
$$\align\l(\f p{Wn+b+h_{s(j)}}\r)=&\da^{(p-1)/2}\l(\f{b+h_{s(j)}}{p}\r)=\da^{(p-1)/2}\l(\f{r_p}p\r)
\\=&\l(\f{r_p\,\delta}p\r)=\cases\da_2&\t{if}\ p=p_{s(i),s(j)},
\\1&\t{otherwise}.\endcases
\endalign$$
Recall that $p_{s(i),s(j)}\|h_{ij}^*$. Therefore,
$$\l(\f{p_{N+i}}{p_{N+j}}\r)=\da\l(\f{h_{ij}^*}{Wn+b+h_{s(j)}}\r)=\da\da_2=\da_1$$
and
$$\l(\f{p_{N+j}}{p_{N+i}}\r)=\da\l(\f{p_{N+i}}{p_{N+j}}\r)=\da_2.$$
This concludes the proof of Theorem 1.3. \qed

\heading{3. Proof of Theorem 1.11}\endheading

The following lemma is a slight modification of Lemma 2.2 which can be proved in a similar way.

\proclaim{Lemma 3.1} Let $k>1$ be an integer. Then there is an admissible set $\Cal H=\{h_1,\ldots,h_k\}$ with $h_1=0<h_2<\ldots<h_k$
such that:

{\rm (i) All those $h_1,h_2,\ldots,h_k$ are multiples of $K=4\prod_{p<4k}p$.

{\rm (ii)} Each $h_i-h_j$ with $1\ls i<j\ls k$ has a prime divisor $p>4k$ with $h_i\not\eq h_j\pmod{p^2}$.

{\rm (iii)} If $1\ls i<j\ls k$, $1\ls s<t\ls k$ and $\{i,j\}\not=\{s,t\}$, then no prime $p>4k$ divides
both $h_i-h_j$ and $h_s-h_t$.
\endproclaim

\proclaim{Lemma 3.2} Let $k>1$ be an integer, and let $\Cal H=\{h_1,\ldots,h_k\}$ with $h_1=0<h_2<\cdots<h_k$ be an admissible set satisfying {\rm (i)-(iii)} in Lemma $3.1$.
Then there is a positive integer $b$ with all of the following properties:

{\rm (i)} $\prod_{i=1}^k(b+h_i)$ is relatively prime to the least common multiple $W$ of those $h_j-h_i$ with $1\ls i<j\ls k$ and $\prod_{2<p\ls w}p$ if $w$ is large enough.

{\rm (ii)} $\prod_{i=1}^k(b+h_i-1)$ is relatively prime to $\prod_{2<p\ls w}p$ if $w$ is large enough.

{\rm (iii)} For any $i,j\in\{1,\ldots,k\}$ with $i\not=j$,  we have $$\l(\f{h_i-h_j}{b+h_j}\r)=-1.$$

{\rm (iv)} If $n>b$, $n\eq b\pmod W$ and $a\in\{h_1,h_1+1,\ldots,h_k\}\sm\Cal H$, then $n+a$ is not prime.
\endproclaim
\Proof. For any $1\ls i <j\ls k$, the number $h_i-h_j$ has a prime divisor $p_{ij}>4k$ with $h_i\not\eq h_j\pmod{p_{ij}^2}$.
Suppose that $p>4k$ is a prime dividing $\prod_{1\ls i<j\ls k}(h_i-h_j)$. Then there is a unique pair $\{i,j\}$ with $1\ls i<j\ls k$
such that $h_i\eq h_j\pmod p$. Note that $p\ls h_k$. As $h_i-h_j\eq h_i-h_i\pmod p$, we have
$$\align&|\{0\ls r<p:\ r\not\eq h_i-h_s,h_i-h_s+1\pmod p\ \t{for all}\ s=1,\ldots,k\}|
\\=&|\{0\ls r<p:\ r\not\eq h_i-h_s,h_i-h_s+1\pmod p\ \t{for all}\ s\in\{1,\ldots,k\}\sm\{i\}\}|
\\\ls& 2(k-1)<\f{p-1}2.
\endalign$$
Recall that $\{1,\ldots,p-1\}$ contains exactly $(p-1)/2$ quadratic residues modulo $p$ and $(p-1)/2$ quadratic nonresidues modulo $p$. So
there is an integer $r_p\not\eq h_i-h_s,\, h_i-h_s+1\pmod p$ for all $s=1,\ldots,k$ such that
$$\l(\f{r_p}p\r)=\cases -(-1)^{\ord_3(h_j-h_i)}&\t{if}\ p=p_{ij},
\\1&\t{otherwise},\endcases$$
where $\ord_3(h_i-h_j)$ is the largest nonnegative integer $a$ such that $3^a|(h_i-h_j)$.
So, for any integer $b\eq r_p-h_i\pmod p$, we have $b+h_s\not\eq0,1\pmod p$.

Assume that
$$S=\{h_1<a<h_k: a\not = h_s, h_s-1\ \t{for all}\ s=1,\ldots,k\}=\{a_i:\ i=1,\ldots,t\}.$$
Clearly $t\ls h_k-k$ and hence we may
choose $t$ distinct primes $q_1,\ldots,q_t\in (h_k,w]$ if $w$ is large enough.
If $b\eq -a_i\pmod {q_i}$, then $b+h_s\eq h_s-a_i\not\eq0,1\pmod {q_i}$ for all $s=1,\ldots,k$ since $|h_s-a_i|<h_k<q_i$.

Let
$$Q=\bg\{p\in(4k,w]:\ p\ \t{is prime and}\ p\nmid\prod_{1\ls i<j\ls k}(h_i-h_j)\bg\}\sm\{q_i:\ i=1,\ldots,t\}.$$
For any prime $q\in Q$, there is an integer $r_q\not\eq-h_i,-h_i+1\pmod q$ for all $i=1,\ldots,k$ since $q>2k$.

By the Chinese Remainder Theorem, there is a positive integer $b$ satisfying the following (1)-(4).

(1) $b\eq 17\pmod {24}$, and $b\eq4\pmod p$ for all primes $p\in[5,4k]$.

(2) $b\eq r_p-h_i\eq r_p-h_j\pmod p$ if $p>4k$ is a prime dividing $h_i-h_j$ with $1\ls i<j\ls k$.

(3) $b\eq -a_i\pmod {q_i}$ for all $i=1,\ldots,t$.

(4) $b\eq r_q\pmod q$ for all $q\in Q$.

By the above analysis,  $\prod_{s=1}^k(b+h_s)(b+h_s-1)$ is relatively prime to $\prod_{2<p\ls w}p$.
Note that $b+h_i\eq 17+0\pmod{24}$ for all $i=1,\ldots,k$. If $w\gs h_k$, then any prime divisor of $W$ does not exceed $w$. So parts (i) and (ii) of Lemma 3.2 are valid.

For each $s=1,\ldots,k$, clearly $b+h_s\eq 17+0\eq1\pmod 8$ and hence
$$\l(\f{-1}{b+h_s}\r)=\l(\f{2}{b+h_s}\r)=1;$$
also $b+h_s\eq 17\eq5\pmod{12}$ and hence
$$\l(\f 3{b+h_s}\r)=\l(\f{b+h_s}3\r)=\l(\f53\r)=-1.$$

Let $i,j\in\{0,\ldots,m\}$ with $i\not=j$. Then
$$\l(\f{h_i-h_j}{b+h_j}\r)=\l(\f{h_{ij}}{b+h_j}\r),$$
where $h_{ij}$ is the odd part of $|h_i-h_j|$. For any prime divisor $p$ of $h_{ij}$, clearly $p\ls h_k\ls w$ and
$$\l(\f p{b+h_j}\r)=\l(\f{b+h_j}p\r).$$
If $3<p<4k$, then $p\mid K$, hence $b+h_j\eq4+0\pmod p$ and thus
$$\l(\f p{b+h_j}\r)=\l(\f{b+h_{j}}{p}\r)=\l(\f{4}p\r)=1.$$
If $p>4k$, then by the choice of $b$ we have
$$\align\l(\f p{b+h_{j}}\r)=&\l(\f{b+h_j}{p}\r)=\l(\f{r_p}p\r)
\\=&\cases-(-1)^{\ord_3|h_j-h_i|}&\t{if}\ p=p_{\min\{i,j\},\max\{i,j\}},
\\1&\t{otherwise}.\endcases
\endalign$$
Recall that $p_{\min{i,j},\max{i,j}}\|h_{ij}$. Therefore,
$$\l(\f{h_i-h_j}{b+h_j}\r)=\l(\f{h_{ij}}{b+h_{j}}\r)=\l(\f 3{b+h_j}\r)^{\ord_3|h_i-h_j|}\l(\f{p_{\min\{i,j\},\max\{i,j\}}}{b+h_j}\r)=-1.$$
So part (iii) of Lemma 3.2 also holds.

Now suppose that $n>b$ is an integer with $n\eq b\pmod W$, and that $a\in\{h_1,h_1+1,\ldots,h_k\}\sm \Cal H$.
If $a=h_s-1$ for some $1\ls s\ls k$, then $n+a\eq b+h_s-1\eq0\pmod 4$ and hence $n+a$ is not prime.
If $a\not= h_s-1$ for all $s=1,\ldots,k$, then $a=a_i$ for some $1\ls i\ls t$, hence
$n+a\eq b+a_i\eq0\pmod {q_i}$ and thus $n+a$ is not prime. (Note that $n+a>W>w\gs q_i$.)
Thus part (iv) of Lemma 3.2 also holds.

In view of the above, we have completed the proof of Lemma 3.2. \qed

\medskip
\noindent{\it Proof of Theorem 1.11}. Choose an integer $k>m$ as in Pollack [P] in the spirit of Maynard-Tao's work.
Let $\Cal H=\{h_1,h_2,\ldots,h_k\}$ be an admissible set constructed in Lemma 3.1 and choose an integer $b$ as in Lemma 3.2.
Let $x$ be sufficiently large, and let $W$ be the least common multiple of those $h_j-h_i$  ($1\ls i<j\ls k$) and $\prod_{2<p\ls \log\log\log x}p$.
Then we have an analogue of Pollack [P, Lemma 3.3].
When $n+h_i$ and $n+h_j$ ($i\not=j$) are both prime with $n\eq b\pmod W$,
we see that $n+h_i$ is a primitive root modulo $n+h_j$ if and only if $|h_i-h_j|$ is a primitive root modulo $n+h_j$
since $n+h_j\eq1\pmod 4$.

Let $P$ be the set of all primes. For $j=1\ldots,k$, set
$$P_j:=\{p\in P:\ |h_i-h_j|\ \t{is a primitive root modulo}\ p \ \t{for every}\ i\not=j\}.$$
Define the weight function $w(n)$ as in [M, Proposition 4.1] or [P, Proposition 3.1], and let $\chi_A(x)$ be the characteristic function of the set $A$. We need to show that
$$\sum\Sb x\ls n\ls 2x\\n\eq b\pmod W\endSb\(\sum_{j=1}^k\chi_{P_j}(n+h_j)\)w(n)\sim \sum\Sb x\ls n\ls 2x\\n\eq b\pmod W\endSb\(\sum_{j=1}^k\chi_{P}(n+h_j)\)w(n).\tag 3.1$$

For a prime $q$ and an integer $g$, define
$$P_q(g)=\l\{p\in P:\ p\eq1\pmod q\ \t{and}\ g^{(p-1)/q}\eq1\pmod p\r\}$$
and
$$\Cal P_q(g)=P_q(g)\sm\bigcup_{q'<q} P_{q'}(g).$$
Under the GRH, Pollack [P, the estimation of $\Sigma_1 - \Sigma_4$] used an effective version of the Chebotarey density theorem to show that if $(\f g{b+h_j})=-1$ then
$$\sum_{q\in P}\sum\Sb x\ls n\ls 2x\\ n\eq b\pmod W\endSb\chi_{\Cal P_q(g)}(n+h_j)w(n)=o\l(\f{\varphi(W)^k}{W^{k+1}}x(\log x)^k\r),$$
where the error term is controlled via the GRH.
Note that if $n\in P\setminus P_j$ then $n\in\Cal P_q(|h_i-h_j|)$ for some $i\neq j$ and some prime $q$. Hence
$$\align&\sum\Sb x\ls n\ls 2x\\n\eq b\pmod W\endSb\(\sum_{j=1}^k\chi_P(n+h_j)\)w(n)
-\sum\Sb x\ls n\ls 2x\\n\eq b\pmod W\endSb\(\sum_{j=1}^k\chi_{P_j}(n+h_j)\)w(n)
\\\ls&\sum_{j=1}^k\sum^k\Sb i=1\\ i\not=j\endSb\sum_{q\in P}\sum\Sb x\ls n\ls 2x\\n\eq b\pmod W\endSb\chi_{\Cal P_q(|h_i-h_j|)}(n+h_j)w(n)
=o\l(\f{\varphi(W)^k}{W^{k+1}}x(\log x)^k\r).
\endalign$$

Maynard and Tao (cf. [M]) have proved
$$
\sum_{\Sb x\ls n\ls 2x\\ n\equiv b\pmod{W}\endSb}w(n)\sim \frac{\alpha\varphi(W)^k}{W^{k+1}}x(\log x)^k\tag 3.2
$$
and
$$
\sum_{\Sb x\ls n\ls 2x\\ n\equiv b\pmod{W}\endSb}\bigg(\sum_{j=1}^k\chi_P(n+h_j)\bigg)w(n)\sim \frac{\beta k\varphi(W)^k}{W^{k+1}}x(\log x)^k,\tag 3.3
$$
where $\alpha$ and $\beta$ are positive constants only depending on $k$ and $w$.
It follows from (3.1) and (3.3) that
$$\sum_{\Sb x\ls n\ls 2x\\ n\equiv b\pmod{W}\endSb}\bigg(\sum_{j=1}^k\chi_{P_j}(n+h_j)\bigg)w(n)\sim \frac{\beta k\varphi(W)^k}{W^{k+1}}x(\log x)^k.$$
Similarly, for each $j=1,\ldots,k$ we have
$$\sum_{\Sb x\ls n\ls 2x\\ n\equiv b\pmod{W}\endSb}\chi_{P_j}(n+h_j)w(n)\sim \frac{\beta \varphi(W)^k}{W^{k+1}}x(\log x)^k,$$
which implies that the set $P_j$ cannot be finite.
Moreover, in view of [M], we may choose a sufficiently large integer $k$ and a suitable weight function $w$ such that
$$\beta k>m\alpha,$$
i.e.,
$$\sum_{\Sb x\ls n\ls 2x\\ n\equiv b\pmod{W}\endSb}\bigg(\sum_{j=1}^k\chi_{P_j}(n+h_j)-m\alpha\bigg)w(n)>0.$$
Since $w(n)$ is non-negative, for some $n\in[x,2x]$ with $n\equiv b\pmod{W}$, $\{n+h_1,\ldots,n+h_k\}$ contains at least $m+1$ primes $n+h_j\ (j\in J)$ with $|J|>m$ and $n+h_j\in P_j$ for all $j\in J$.
According to the construction of $b$ and Lemma 3.2 (iv), for each $j=1,\ldots,k$, the interval $(n+h_j,n+h_{j+1})$ contains no prime. So those primes in $\{n+h_1,\ldots,n+h_k\}$ are consecutive primes. Note that
$n+h_k-(n+h_1)=h_k-h_1$ is a constant $C_m$ only depending on $m$.
For any $i,j\in J$ with $i\not=j$, the number $h_i-h_j$, as well as the prime $n+h_i$, is a primitive root modulo the prime $n+h_j$. This concludes the proof. \qed

\Ack. The authors would like to thank the referee for helpful comments.


\widestnumber\key{BFTB}

 \Refs


\ref\key BFTB\by W.D. Banks, T. Freiberg and C.L. Turnage-Butterbaugh\paper Consecutive primes in tuples\jour Acta Arith.\vol 167\yr 2015\pages 261--266\endref

\ref\key M\by J. Maynard\paper Small gaps between primes\jour Annals of Math.\vol 181\yr 2015\pages 383--413\endref

\ref\key P\by P. Pollack\paper Bounded gaps between primes with a given primitive root\jour Algebra Number Theory\vol 8\yr 2014\pages 1769--1786\endref

\ref\key Po\by D.H.J. Polymath\paper Bounded gaps between primes\jour Polymath 8 Project. Available from
{\tt http://michaelnielsen.org/polymath1/index.php?title=Bounded\_gaps\_between\_primes}\endref

\ref\key S\by D.K.L. Shiu\paper Strings of congruent primes\jour J. London Math. Soc.\vol 61\yr2000\pages 359--373\endref

\ref\key Su\by Z.-W. Sun\paper {\rm Sequences A243837, A243839 and A243901 in OEIS (On-Line Encyclopedia of Integer Sequences)}
\jour {\tt http://oeis.org}\endref

\ref\key Z\by Y. Zhang\paper Bounded gaps between primes\jour Annals of Math.\vol 179\yr 2014\pages 1121--1174\endref

\endRefs
\enddocument